\documentclass[11pt]{amsart}

\usepackage{amsmath,amsthm, amscd, amssymb, amsfonts}
\usepackage[all]{xy}
\usepackage{verbatim}

\newtheorem{theorem}{Theorem}[section]
\newtheorem*{maintheorem}{Main Theorem}

\newtheorem{lem}[theorem]{Lemma}
\newtheorem{cor}[theorem]{Corollary}

\newtheorem{pro}[theorem]{Proposition}

\theoremstyle{definition}
\newtheorem{fed}[theorem]{Definition}

\newtheorem{exas}[theorem]{Examples}

\theoremstyle{remark}
\newtheorem{rem}[theorem]{Remark}
\newtheorem{rmk}[theorem]{Remarks}
\newtheorem*{proofmt}{Proof of the Main Theorem}

\newcommand{\ydh}{{}^{H}_{H}\mathcal{YD}}
\newcommand{\ydg}{{}^{G}_{G}\mathcal{YD}}
\newcommand{\ydho}{{}^{H_0}_{H_0}\mathcal{YD}}

\newcommand\id{\operatorname{id}}

\newcommand\End{\operatorname{End}}
\newcommand\sg{\operatorname{sgn}}

\newcommand\aut{\operatorname{Aut}}

\newcommand\Gr{\operatorname{Gr}}
\newcommand\gr{\operatorname{gr}}
\newcommand\co{\operatorname{co}}
\newcommand\ad{\operatorname{ad}}

\def\k{\Bbbk}

\def\rhdi{\rhd^{-1}}
\def\ot{\otimes}

\def\s{\mathbb{S}}

\def\B{\mathfrak{B}}

\def\eps{\epsilon}

\def\mD{\mathcal{D}}
\def\mB{\mathcal{B}}

\def\mO{\mathcal{O}}
\def\mR{\mathcal{R}}
\def\mQ{\mathcal{Q}}
\def\mH{\mathcal{H}}
\def\hq{\mH(\mQ)}

\newcommand{\J}{{\mathcal J}}

\def\pf{\begin{proof}}
\def\epf{\end{proof}}

\begin{document}

%%%%%%%%% TITLE  AND OTHERS %%%%%%%%%%%%%%%%%%%%%

\title[Finite dimensional pointed Hopf algebras over $\s_4$]
{Finite dimensional pointed Hopf algebras \\ over $\s_4$}

\author[Garc\'ia, Garc\'ia Iglesias]{Gast\'on Andr\'es Garc\'ia; \, Agust\'in Garc\'ia Iglesias}

\address{FaMAF-CIEM (CONICET), Universidad Nacional de C\'ordoba,
Medina A\-llen\-de s/n, Ciudad Universitaria, 5000 C\' ordoba, Rep\'
ublica Argentina.} \email{(ggarcia|aigarcia)@mate.uncor.edu}

\thanks{The work was partially supported by CONICET,
FONCyT-ANPCyT, Secyt (UNC)}

\date{\today}

\begin{abstract}
Let $\k$ be an algebraically closed field of characteristic $0$.
We conclude the classification of finite-dimensional pointed Hopf
algebras whose group of group-likes is $\s_4$. We also describe
all pointed Hopf algebras over $\s_5$ whose infinitesimal braiding
is associated to the rack of transpositions.
\end{abstract}

\maketitle

%%%%%%%%%%%%%%%%%%%%%%%%%%%%%%%%%%%%%%%%%%%%%%%%%%%%%%%%%%%

%%%%%%%%% PAPER %%%%%%%%%%%%%%%%%%%%%%%%%%%%%%%%%%%%%%%%%%%

%%%%%%%%%%%% INTRODUCTION Y PRELIMINARES %%%%%%%%%%%%%%%%%%
\section{Introduction}

The classification of finite-dimensional Hopf algebras over an
algebraically closed field $\k$ of characteristic $0$ is a widely
open problem. Strictly different techniques are employed when
dealing with semisimple or non semisimple algebras. In this later
case, useful techniques have been developed for the special case
of pointed Hopf algebras.

\smallbreak A significant progress has been achieved in \cite{AS2}
in the case of pointed Hopf algebras with abelian group of
group-likes. When the group of group-likes is not abelian, the
problem is far from being completed. Some hope is present in the
lack of examples: in this situation, Nichols algebras tend to be
infinite dimensional, see for example \cite{AF}, \cite{AFGV},
\cite{AZ}, \cite{AFZ}. Nevertheless, examples on which the Nichols
algebras are finite dimensional do exist. Over $\s_4$ these
algebras are determined in \cite{AHS}: there are 3 of them, all
arising from racks in $\s_4$, associated to a cocycle. When the
cocycle is $-1$, there are two Nichols algebras, one corresponding
to the rack of transpositions and the other to the rack of
4-cycles, presented in \cite{MS} and \cite{AG2}, respectively.
When the cocycle is non-constant, this Nichols algebra was defined
in \cite{MS} and, independently, in \cite{FK} as a quadratic
algebra. These three Nichols algebras are an exhaustive list of
Nichols algebras in the category of Yetter Drinfeld modules
$_{\s_4}^{\s_4}\mathcal{YD}$, and in this paper we prove that all
pointed Hopf algebras $H$ with $G(H)\cong\s_4$ are in fact
liftings of them. Up to now, only two Nichols algebras over $\s_5$
are known to be finite-dimensional. Here, we describe all their
liftings.

\smallbreak In \cite{AG1}, a family of pointed Hopf algebras was
defined; these were shown, by means of Gr\"obner basis, to be
liftings of the Nichols algebras associated to the class of
transpositions and constant cocycle -1. Following ideas from that
paper, we define two more families of Hopf algebras, and show that
they are liftings of the other Nichols algebras. We remark that
our proofs do not use Gr\"obner basis, instead, we develop part of
the representation theory for these algebras and combine this new
knowledge with some techniques on quadratic algebras to prove that
these algebras are liftings of $\s_4$, see Proposition
\ref{pro:dimension}. The main theorem of the present paper is the
following:
\begin{maintheorem}\label{main}
Let $H$ be a finite-dimensional pointed Hopf algebra with
$G(H)\cong\s_4$, $H\neq\k\s_4$. Then $H$ is isomorphic to one and
only one of the following algebras:
\begin{enumerate}
\item $\B(\mO_2^4,-1)\sharp\k\s_4$;
\item $\B(\mO_4^4,-1)\sharp\k\s_4$;
\item $\B(\mO_2^4,\chi)\sharp\k\s_4$;
\item $\mH(\mQ_4^{-1}[t])$, for exactly one $t\in\mathbb{P}^1_\k$;
\item $\mH(\mD[t])$, for exactly one $t\in\mathbb{P}^1_\k$;
\item $\mH(\mQ_4^{\chi}[1])$.
\end{enumerate}
\end{maintheorem}
See Subsection 2.4 and Definitions \ref{aene}--\ref{ace} for a
description of the algebras mentioned above. We point out that
$\s_4$ is the second finite non-abelian group which admits
non-trivial finite-dimensional pointed Hopf algebras such that the
classification is completed, the first one being $\s_3$ \cite{AHS}.

The paper is organized as follows. In Section 2 we recall some basic
facts about Hopf algebras and racks, with examples, together with
the definitions of three  Nichols algebras in
$_{\s_4}^{\s_4}\mathcal{YD}$ and two in
$_{\s_5}^{\s_5}\mathcal{YD}$. In fact, we present a proof showing
that two of these braided Hopf algebras are Nichols algebras, a fact
that seems to be absent in the literature. All of these algebras are
quadratic. We explicitly describe the set of quadratic relations in
a Nichols algebra associated to a rack and a cocycle. In Section 3
we show that a pointed Hopf algebra whose infinitesimal braiding is
subject to certain hypothesis is generated by group-likes and
skew-primitives. We define families of pointed Hopf algebras
associated to quadratic lifting data (Def. \ref{fed:ql-datum}) and
explicitly distinguish three of these families, which will be of our
interest later on. In Section 4 we use some quadratic-algebra
techniques to provide a bound on the dimension of an algebra
belonging to one of the families defined previously, and we use this
to prove that, under certain conditions, these families give all
possible pointed Hopf algebras over a given group. In Section 5 we
develop part of the representation theory for the algebras we have
remarked in Section 3, in order to prove that for each one the group
of group-likes is the symmetric group $\s_n$. Thus we show that
these algebras are liftings of the Nichols algebras over $\s_n$ from
Section 2. Finally, in Section 6 we prove the Main Theorem. We also
show that the liftings of $\s_5$ defined exhaust the list of
finite-dimensional pointed Hopf algebras over $\s_5$ which are
liftings of one of the known finite-dimensional Nichols algebras
over this group.

\section{Preliminaries}
\subsection{Conventions}Let $H$ be a Hopf
algebra over $\k$, with antipode $S$. We will use the Sweedler's
notation $\Delta(h)=h_1\ot h_2$ for the comultiplication \cite{S}.
Let $\ydh$ be the category of (left-left) Yetter-Drinfeld modules
over $H$. That is, $M$ is an object of $\ydh$ if and only if there
exists an action $\cdot$ such that $(M,\cdot)$ is a (left)
$H$-module and a coaction $\delta$ such that $(M,\delta)$ is a
(left) $H$-comodule, subject to the following compatibility
condition:
\begin{equation*}
\delta(h\cdot m)=h_1m_{-1}S(h_3)\ot h_2\cdot m_{0}, \ \forall\, m\in
M, h\in H,
\end{equation*}
where $\delta(m)=m_{-1}\ot m_0$. If $G$ is a finite group and
$H=\k G$, we denote $\ydg$ instead of $\ydh$.

\subsection{The Lifting Method for the classification of Hopf algebras}

Let $(V,c)$ be a braided vector space, \emph{i.~e.} $V$ is a
vector space and $c\in GL(V\otimes V)$ is a \emph{braiding}:
$(c\ot \id)(\id\ot c)(c\ot \id) = (\id\ot c)(c\ot \id)(\id\ot
c)\in\End(V\ot V\ot V)$. Let $\B(V)$ be the Nichols algebra of
$(V,c)$, see \emph{e.~g.} \cite{AS1}. Let $\J=\J_V$ be the kernel
of the canonical projection $T(V) \to \B(V)$; $\J = \oplus_{n\geq
0}\J^n$ is a graded Hopf ideal. If $\mathcal{J}_2$ is the ideal
generated by $\J^2$, $\widehat\B_2(V)=T(V)/\J_2$ is the
\emph{quadratic} Nichols algebra of $(V,c)$.

The \emph{coradical} $H_0$ of $H$ is the sum of all simple
sub-coalgebras of $H$. In particular, if $G(H)$ denotes the group
of \emph{group-like elements} of $H$, we have $\k G(H)\subseteq
H_0$. We say that a Hopf algebra is \emph{pointed} if $H_0=\k
G(H)$. Define, inductively, $H_n=\Delta^{-1}(H\ot H_{n-1}+H_0\ot
H)$. The family $\{H_i\}_{i\geq 0}$ is the \emph{coradical
filtration} of $H$, it is a coalgebra filtration and satisfies
$H_i\subseteq H_{i+1}$, $H=\bigcup H_i$, see \cite{M}, \cite{S}.
It is a Hopf algebra filtration if and only if $H_0$ is a sub-Hopf
algebra of $H$. In this case, for $\gr H(n) = H_n/H_{n-1}$ (set
$H_{-1} =0$), $\gr H = \oplus_{n\ge 0}\gr H(n)$ is the associated
graded Hopf algebra. Let $\pi:\gr H\to H_0$ be the homogeneous
projection. We have the following invariants of $H$:

\begin{itemize}
 \item $R= (\gr H)^{\co \pi}$ is the
\emph{diagram} of $H$; this is a braided Hopf algebra in $\ydho$,
and it is a graded sub-object of $\gr H$.

\item $V:=R(1)=P(R)$, with the braiding from $\ydho$, is called the \emph{infinitesimal braiding} of $H$.

\end{itemize}

It follows that the Hopf algebra $\gr H$ is the Radford biproduct
$\gr H \simeq R\# \k G(H)$. The subalgebra of $R$ generated by $V$
is isomorphic to the Nichols algebra $\B(V)$.

Let $\Gamma$ be a finite group and let $H_0$ be the group algebra of
$\Gamma$. The main steps of the \emph{Lifting Method} \cite{AS1} for
the classification of all finite-dimensional pointed Hopf algebras
with group $\Gamma$ are:

\begin{itemize}
\item determine all $V\in \ydho$ such that the Nichols algebra $\B(V)$ is finite dimensional,

\item for such $V$, compute all Hopf algebras $H$ such
that $\gr H\simeq \B(V)\sharp \k\Gamma$. We call $H$ a
\emph{lifting} of $\B(V)$ over $\Gamma\cong G(H)$.

\item Prove that any finite-dimensional pointed Hopf algebras with group $\Gamma$
is generated by group-likes and skew-primitives.
\end{itemize}

\subsection{Racks} A \emph{rack} is a pair $(X,\rhd)$, where $X$ is a
non-empty set and $\rhd:X\times X\to X$ is a function, such that
$\phi_i=i\rhd (\cdot):X\to X$ is a bijection for all $i\in X$
satisfying:
$$
i\rhd(j\rhd k)=(i\rhd j)\rhd (i\rhd k), \, \forall i,j,k\in X.
$$
The following are examples of racks.
\begin{itemize}
\item Let $G$ be a group, $g\in G$, and $\mO_g$ the conjugacy
class of $g$. Then, if $x\rhd y=xyx^{-1}$, $(\mO_g, \rhd)$ is a
rack.
\item If $G=\s_n$ and $\tau$ is a $j$-cycle, we will denote by
$\mO_j^n$, or $\mO_j$ if $n$ is fixed by the context, the rack
induced by its conjugacy class.
\item In particular, we will work with the racks $\mO_2^n$ of
transpositions in $\s_n$ and with the rack $\mO_4^4$ of $4$-cycles
in $\s_4$.
\item If $(X,\rhd)$ is a rack, then the inverse rack $(X^{-1},\rhd^{-1})$, see \cite{afgv2}, is
given by $X^{-1}=X$ and $i\rhd^{-1}j=k$, provided $i\rhd k=j$,
$\forall\,i,j,k\in X^{-1}$. Indeed, note that $i\rhdi (j\rhdi
k)=s\Leftrightarrow i\rhd s=j\rhdi k \Leftrightarrow j\rhd (i\rhd
s)=k$. While $(i\rhdi j)\rhdi (i\rhdi k)=u\Leftrightarrow (i\rhdi
j)\rhd u=i\rhdi k\Leftrightarrow i\rhd((i\rhdi j)\rhd
u)=k\Leftrightarrow(i\rhd(i\rhdi j))\rhd (i\rhd u)=k\Leftrightarrow
j\rhd (i\rhd u)=k\Leftrightarrow i\rhd u=i\rhd s\Leftrightarrow
u=s$, since $\phi_j,\phi_i$ are bijections.
\end{itemize}

Note that the racks $\mO_4^4$ and $\mO_2^4$ on $\s_4$ are not
isomorphic, since, for $i\in\mO_2^4$, we always have $\phi_i^2=\id$,
and for $i\in\mO_4^4$ $\phi_i^4=\id$, but $\phi_i^2\neq\id$.

A rack $(X,\rhd)$ is said to be \emph{indecomposable} if it cannot
be decomposed as the disjoint union of two sub-racks. It is said
to be \emph{faithful} if $\phi_i=\phi_j$ only for $i=j$. We refer
the reader to \cite{AG2} for detailed information on racks.

Let $(X,\rhd)$ be a rack. A 2-cocycle $q:X\times X\to \k^*$,
$(i,j)\mapsto q_{ij}$ is a function such that
\begin{equation*}\label{cond cociclo} q_{i,j\rhd k}q_{j,k}=q_{i\rhd j,i\rhd
k}q_{i,k}, \ \forall\,i,j,k\in X.
\end{equation*}
In this case it is possible to generate a braiding $c^q$ in the
vector space $\k X$ with basis $\{x_i\}_{i\in X}$ by $c^q(x_i\otimes
x_j)=q_{ij}x_{i\rhd j}\otimes x_i, \, \forall\,i,j\in X.$

\begin{exas}\label{q}\label{chi}
\begin{itemize}
\item Let $(X,\rhd)$ be a rack. A constant map $q:X\times
X\to \k^*$, \emph{i.~e.} $q_{\sigma\tau}=\xi\in\k^*$ for every
$\sigma,\tau\in X$, is a $2$-cocycle.

\medbreak

\item\cite[Ex. 5.3]{MS} Let $(X,\rhd)=\mO_2^n$. A
$2$-cocycle is given by a function $\chi:X\times X\to \k^*$ defined
as, if $\tau, \sigma\in X$, $\tau=(ij)$ and $i<j$:
\begin{equation*}\label{cociclo} \chi(\sigma, \tau) =
\begin{cases}
  1,  & \mbox{if} \ \sigma(i)<\sigma(j) \\
  -1, & \mbox{if} \ \sigma(i)>\sigma(j).
\end{cases}
\end{equation*}
\end{itemize}
\end{exas}

The dual braided vector space of $(\k X,c^q)$ is given by $(\k
X^{-1},c^{\tilde{q}})$, with
\begin{equation}\label{eqn:qtilde}
\tilde{q}_{kl}=q_{k\, k\rhdi l},\quad k,l\in X^{-1}.
\end{equation}
Indeed, let
$Y=\{y^i : i\in X\}$ be a dual basis of $\{x_i : i\in X\}$, and
consider the pairing $\langle\cdot,\cdot\rangle$ defined as
$\langle x_i\ot x_j, y^k\ot y^l\rangle=\delta_{i,l}\delta_{j,k}$,
for $i,j,k,l\in X$. The braiding $c^*$ is thus $\langle x_i\ot
x_j, c^*(y^k\ot y^l)\rangle =\langle c(x_i\ot x_j), y^k\ot
y^l\rangle= q_{ij}\langle x_{i\rhd j}\ot x_i, y^k\ot
y^l\rangle=q_{ij}\delta_{l,i\rhd j}\delta_{k,i}=q_{k,k\rhdi
l}\delta_{k\rhdi l,j}\delta_{k,i}$ and thus $c^*(y^k\ot
y^l)=q_{k,k\rhdi l}y^{k\rhdi l}\ot y^k$.

\subsection{Nichols algebras associated to racks}\label{subsect:R}
Let $X$ be a rack, $q$ a 2-cocycle. We will denote by $\B(X,q)$
the Nichols algebra associated to the braided vector space $(\k
X,c^q)$. Let $\mR$ be the set of equivalence classes in $X\times
X$ for the relation generated by $(i,j)\sim(i\rhd j,i)$. Let
$C\in\mR$, $(i,j)\in C$. Take $i_1=j$, $i_2=i$, and recursively,
$i_{h+2}=i_{h+1}\rhd i_h$. For each $C$, there exists
$n\in\mathbb{N}$ such that $i_{n+k}=i_k$, $k\in\mathbb{N}$. Let
$n(C)$ be the minimum such $n$; thus $C=\{i_1, \ldots, i_{n(c)}\}$
and $n(C)=\#C$. Let $\mR'$ be the set of all $C\in \mR$ that
satisfy
\begin{align}\label{eqn:condcocycle}
\prod_{h=1}^{n(C)} q_{i_{h+1},i_h} &=(-1)^{n(C)}.
\end{align}

\begin{lem}\label{lem:quadraticrels}
A basis of the space $\J^2$ of quadratic relations of $\B(X,q)$ is
given by
\begin{align}
\label{eqn:quadraticrelations}   b_C:&= \sum_{h=1}^{n(C)}\eta_h(C)
\, x_{i_{h+1}}x_{i_h}, \quad C\in \mR',
\end{align}
where $\eta_1(C)=1$ and
$\eta_h(C)=(-1)^{h+1}q_{{i_2i_1}}q_{{i_3i_2}}\ldots
  q_{{i_hi_{h-1}}}$, $h\ge 2$.
\end{lem}

\pf Let $U_C$ be the subspace of $\k X \otimes \k X$ spanned by
$x_i\otimes x_j$, $(i,j)\in C$. Then $\J^2 = \ker (c^q+\id) =
\oplus_{C\in \mR} \ker (c^q+\id)\vert_{U_C}$. We have $\det
(c^q+\id)\vert_{U_C} = \prod_{h=1}^{n(C)} q_{i_{h+1},i_h}
(-1)^{n(C) + 1} + 1$. If this is 0, then $b_C$ spans  $\ker
(c^q+\id)\vert_{U_C}$. \epf

\begin{rmk} The constant cocycle $q\equiv-1$ evidently fulfills
\eqref{eqn:condcocycle}, and $\eta_h(C)=1$ for every $C\in\mR,
h=1, \ldots,   n(C)$.

\medbreak More generally, the constant cocycle $q\equiv\omega$, for
$-\omega$ a primitive $l$-th root of unity and $l|n(C)$, also
fulfills \eqref{eqn:condcocycle}, and $\eta_h(C)=(-\omega)^{h-1}$
for every $C\in\mR, h=1, \ldots,   n(C)$. See \cite[Lem. 6.13]{AG2}.

\medbreak It is easy to see that the cocycle $\chi$ in Examples
\ref{chi} fulfills \eqref{eqn:condcocycle} for every $C\in\mR$.

\medbreak When dealing with the racks $\mO_2^n,
  \mO_4^4$, an equivalence class $C\in\mR$ may have only 1, 2, or 3
  elements, and thus a relation on the basis is composed of at most
  three summands.
\end{rmk}

\begin{theorem}\label{nichols:cte}
\begin{enumerate}
\item Let $3\leq n\leq 5$. The Nichols algebras
$\B(\mO_2^n,-1)$ are finite dimensional, of dimensions 12, 576,
8294400, respectively, and it holds
$\B(\mO_2^n,-1)=\widehat\B_2(\mO_2^n, -1)$.
\item $\dim\B(\mO_4^4,-1)=576$ and
$\B(\mO_4^4,-1)=\widehat\B_2(\mO_4^4, -1)$.
\end{enumerate}
\end{theorem}
\begin{proof}
We see (1). Cases $n=3,4$ are in \cite[Ex. 6.4]{MS}. For $n=5$, the
algebra was introduced as a quadratic algebra in \cite{MS}. Gra\~na
established $\widehat\B_2(V)=\B(V)$ in \cite{G2}. Its dimension was
determined by Roos with computer program \verb"Bergman". (2) is
\cite[Th. 6.12]{AG2}, using \verb"Bergman".
\end{proof}

It turns out that the Nichols algebras $\B(\mO_2^4, \chi)$,
$\B(\mO_2^5, \chi)$ are also quadratic. To see this, we recall here
the notion of differential operators associated to a rack $(X,\rhd)$
in a Nichols algebra $\B(X)$ arising from that rack. For $x\in X$,
let $x^*$ be the dual element to $x$ in $(\k X)^*$ and, for $m\neq
1$, extend it as zero in $\B^m$, the homogeneous component of degree
$m$ of $\B(V)$. Define $\delta_x:\B^n\to\B^{n-1}$ by
$\delta_x=(\id\ot x^*)\Delta$, where $\Delta$ is the
comultiplication of the Nichols algebra. These skew-derivations
provide a powerful tool to decide whether an element in the algebra
is zero or different from zero, since $\alpha\in\B^m$ $(m\geq 2)$ is
zero if and only if $\delta_x(\alpha)=0$ $\forall\,x\in \k X$, see
\cite[Section 6]{AG2}. The following proposition has been assumed to
be true in the literature. We thank Mat\'ias Gra\~na for providing
us the computations needed to finish its proof.

\begin{pro}\label{nichols:milsch}
Let $n=4,5$. The Nichols algebra $\B(\mO_2^n,\chi)$ is quadratic
and has dimension 576, 8294400, respectively.
\end{pro}
\begin{proof}
Let $n=4$, $B=\widehat\B_2(\mO_2^4,\chi)$. According to
\cite[Probl. 2.3]{FK}, the Hilbert polynomial of the algebra $B$
is
\begin{equation}\label{hp}
P_B(t)=[2]^2[3]^2[4]^2,
\end{equation}
where we denote by $[k]$ the polynomial $1+t+t^2+\ldots+t^{k-1}$.
Thus, $\dim B=P_B(1)=2^23^34^2=576$. Equation \eqref{hp} also
implies that the top degree (cf. \cite[Section 6]{AG2}) of $B$ is
12. If $\B^{12}(\mO_2^4,\chi)\neq 0$ we would get
$B=\B(\mO_2^4,\chi)$ by \cite[Th. 6.4 (2)]{AG2}. Let $$
a=x_{(12)}, \ b=x_{(13)}, \ c=x_{(14)}, \ d=x_{(23)}, \
e=x_{(24)}, \ f=x_{(34)}.
$$Consider $abacabacdedf\in\B^{12}(\k\mO_2^4,\chi)$. We use derivations and the help the computer program \verb"Deriva" \cite{G1} developed by Mat\'ias Gra\~ na to see that
$$
\delta_c\delta_b\delta_c\delta_a\delta_b\delta_d\delta_c\delta_b\delta_f\delta_d\delta_f\delta_e(abacabacdedf)\neq 0
$$
and thus $\B^{12}(\mO_2^4,\chi)\neq 0$.

For $n=5$, the result follows similarly. In this case, the Hilbert
polynomial is $P(t)=[4]^4[5]^2[6]^4$, see again \cite[Probl.
2.3]{FK}. Therefore we have that $\dim
\widehat\B_2(\mO_2^5,\chi)=8294400$ and that its top degree is
$40$. Now, take $a, b, c, d, e, f$ as before and
$$ g=x_{(15)}, \ h=x_{(25)}, \
k=x_{(35)}, \ m=x_{(45)}.
$$
Again, using \verb"Deriva", we see that $\B^{40}(\mO_2^5,\chi)\neq
0$ since the operator
\begin{align*}
\nabla=&\delta_k\delta_m\delta_f\delta_m\delta_e\delta_h\delta_f\delta_e\delta_c\delta_h\delta_e\delta_h\delta_g\delta_k\delta_c\delta_b\delta_m\delta_d\delta_c\\
&\circ\delta_e\delta_m\delta_g\delta_e\delta_k\delta_g\delta_c\delta_h\delta_e\delta_g\delta_k\delta_d\delta_e\delta_g\delta_k\delta_b\delta_a\delta_d\delta_c\delta_a\delta_b
\end{align*}
satisfies $\nabla(abadabadgabadabadgabagdabadgcechcechfkfm)\neq
0$.
\end{proof}
Let us describe explicitly all these Nichols algebras:
\begin{align*}
\B(\mO_2^n,-1)=  \,\k\langle &x_{(ij)}, 1\le i < j \le n \vert
x_{(ij)}^2,\, x_{(ij)}x_{(kl)}+x_{(kl)} x_{(ij)}, \, (ij) \neq (kl), \\
&x_{(ij)}x_{(jk)}+x_{(jk)} x_{(ik)}+x_{(ik)} x_{(ij)}, \quad 1\le
i < j < k \le n \rangle;
\\
\B(\mO_2^n,\chi)=  \,\k\langle &x_{(ij)}, 1\le i < j \le n \vert
x_{(ij)}^2,\, x_{(ij)}x_{(kl)} - x_{(kl)} x_{(ij)}, \, (ij) \neq (kl), \\
&x_{(ij)}x_{(jk)} - x_{(jk)} x_{(ik)} - x_{(ik)} x_{(ij)},
\\&x_{(jk)}x_{(ij)} -  x_{(ik)}x_{(jk)} - x_{(ij)}x_{(ik)}, \quad
1\le i < j < k \le n \rangle;
\\
\B(\mO_4^4,-1) =\, \k\langle &x_\sigma, \sigma\in \mO_4^4 \vert
x_\sigma^2,\, x_\sigma
x_{\sigma^{-1}}+x_{\sigma^{-1}} x_\sigma, \\
& x_\sigma x_\tau+x_\nu x_\sigma+x_\tau x_\nu,  \mbox{if }
\sigma\tau=\nu\sigma \ \mbox{and} \ \tau\neq \sigma\neq \nu\in
\mO_4^4 \rangle.
\end{align*}

\subsection{Nichols algebras of Yetter-Drinfeld modules}
Let $G$ be a finite group, $\mO=\{t_1, \ldots, t_n\}$ a conjugacy
class in $G$ and $\rho: G^s\to\mathrm{GL}(V)$ an irreducible
representation of the centralizer $G^s$ of a fixed element
$s\in\mO$. Assume $s=t_1$, and let $g_i\in G$, $i=1, \ldots, n$ such
that $g_isg_i^{-1}=t_i$. Irreducible Yetter Drinfeld modules over
$G$ are in bijective correspondence with pairs $(\mO, \rho)$. For
such pair, the module $M(\mO,\rho)$ is defined as
$\bigoplus\limits_{i\in I}g_i\ot V$, with action and coaction given
by:
$$ g\cdot (g_i\ot v)=g_j\ot(\rho(\gamma)(v)), \ \delta(g_i\ot
v)=t_i\ot(g_i\ot v),
$$
for $\gamma\in G^s$ such that $gg_i=g_j\gamma$.  When the
representation $\rho$ is one dimensional, the underlying braided
vector space of $M(\mO,\rho)$ is $(\k X,c^q)$, where $X$ is the rack
given by the conjugation on $\mO$ and $q:X\times X\to\k^*$ is a
2-cocycle. The equivalence is given by $q(t_i,t_j)=\rho(\gamma)$,
for $\gamma\in G^s$ such that $g_ig_j=g_k\gamma$.

Let $\B(\mO,\rho)$ denote the Nichols algebra of $M(\mO,\rho)$. The
following theorem lists all irreducible modules $M(\mO,\rho)$ in
$_{\s_4}^{\s_4}\mathcal{YD}$ such that $\B(\mO,\rho)$ is
finite-dimensional. If $G=\s_4, s=(1234)$, let $\rho_\_$ be the
character of $G^s=\langle s \rangle\cong\mathbb{Z}_4$ given by
$\rho_\_(s)=-1$. If $s=(12)$, $G^s=\langle (12), (34)
\rangle\cong\mathbb{Z}_2\times\mathbb{Z}_2$; let $\eps$ and $\sg$ be
the trivial and the sign representations of $\mathbb{Z}_2$,
respectively. Note that $M(\mO_2^n,\sg\ot\sg)\cong (\k\mO_2^n,c^q)$,
$M(\mO_2^n,\sg\ot\eps)\cong (\k\mO_2^n,c^\chi)$ and
$M(\mO_4^4,\rho_\_)\cong (\k\mO_4^4,c^q)$, for $n=4,5$ and
$q\equiv-1$.
\begin{theorem}\cite[Th. 4.7]{AHS}\label{teo:ahs}
The only Nichols algebras of Yetter-Drinfeld modules over $\s_4$
with finite dimension are $\B(\mO_2^4,\sg\ot\sg)$,
$\B(\mO_2^4,\sg\ot\eps)$ and $\B(\mO_4^4,\rho_\_)$, up to
isomorphism. All have dimension 576.\qed
\end{theorem}
\begin{rem}\label{rem:O23}
When $n=5$, if $M\neq M(\mO_2^5,\sg\ot\mu),
M(\mO_{2,3}^5,\sg\ot\eps)$ for $\mu\in\{\eps,\sg\}$, the associated
Nichols algebra is infinite dimensional, see \cite{AFZ}. It is not
known if the Nichols algebra $\B(\mO_{2,3}^5,\sg\ot\eps)$ has finite
dimension; this is the only open case in $\s_5$.
\end{rem}

If $H$ is a Hopf algebra with bijective antipode (in particular, if
$H=\k G$), the category $_H^H\mathcal{YD}$ is rigid. Explicitly, let
$M\in\ydh$ with linear basis $\{e_i\}_{i=1}^n$ and dual basis
$\{e^i\}_{i=1}^n$. Its (right) dual $M^*$ is the dual $\k$-vector
space of $M$ with action and coaction given by:
\begin{equation}\label{dual}
(h\cdot f)(m)=f(S(h)m), \quad f_{(-1)}\ot f_{(0)}=\sum_{i=1}^n
S^{-1}((e_i)_{(-1)})\ot f((e_i)_{(0)})e^i,
\end{equation}
for  $f\in M^*$, $m\in M$. The left dual $^*M$ is defined
analogously. When $H=\k G$, $^*M\cong M^*$. If $X$ is a finite
rack and $q$ a 2-cocycle, $\B^n(X,q)^*\cong
\B^n(X^{-1},\tilde{q})$, for $\tilde{q}$ as in \eqref{eqn:qtilde}.
Besides, if $\dim\B(X,q)<\infty$, then $\B(X,q)^*\cong
\B(X^{-1},\tilde{q})$, see \cite[Lem. 2.6]{G3}.

\section{Pointed Hopf algebras associated to racks.}
Let $X$ be a rack, $q$ a 2-cocycle. In this section we are
interested in  finite-dimensional pointed Hopf algebras $H$ such
that its infinitesimal braiding arises from a \emph{principal
YD-realization} of $(X,q)$ over a finite group $G$. We recall from
\cite[Def. 3.2]{AG1}: this is a collection $(\cdot, g,
(\chi_i)_{i\in X})$ where $\cdot$ is an action of $G$ on $X$;
$g:X\to G$ is a function such that $g_{h\cdot i} = hg_ {i}h^{-1}$
and $g_{i}\cdot j=i\rhd j$; and the family $(\chi_i)_{i\in X}$, with
$\chi_i:G\to\k^*$, is a 1-cocycle, \emph{i.~e.}
$\chi_i(ht)=\chi_i(t)\chi_{t\cdot i}(h)$, for all $i\in X$, $h,t\in
G$, satisfying $\chi_i(g_{j})=q_{ji}$. The realization is said to be
\emph{faithful} if $g$ is injective. Let $K$ be the subgroup of $G$
generated by the image of $X$; $K$ acts by rack automorphisms on
$X$. In Theorem \ref{conj} below we show that, if $\B(X,q)$ is
finite-dimensional and $\B(X^{-1},\tilde{q})$ is quadratic, such a
pointed Hopf algebra must be generated by group-likes and skew
primitive elements. Another relevant result concerning these
algebras is Lemma \ref{lem:lambda}, where we show that certain
relations hold in such an $H$. Finally, we define families of
pointed Hopf algebras associated to a given data, some of which will
be shown to be liftings of quadratic Nichols algebras.

\subsection{Generation in degree one}
The following theorem, which agrees with a well-known conjecture
\cite[Conj. 1.4]{AS}, is a key step in the classification. When
$(X,q)=(\mO_2^n,-1)$, $n=3,4,5$ the result is \cite[Th. 2.1]{AG1}.

\begin{theorem}\label{conj}
Let $(X,q)$ be a rack such that $\B(X,q)$ is finite-dimensional and
$\B(X^{-1},\tilde{q})$ is quadratic. Let $H$ be a
pointed Hopf algebra such that its infinitesimal braiding arises
from a principal YD-realization of $(X,q)$ over a finite group $G$.
Then the diagram of $H$ is a Nichols algebra, and consequently $H$
is generated by its group-like and skew-primitive elements.
\end{theorem}
\begin{proof}
Let $G=G(H)$, $\gr H=R\sharp\k G$, $S=R^*$ the graded dual. By
\cite[Lem. 5.5]{AS}, $S$ is generated by $V=S(1)$ and $R$ is a
Nichols algebra if and only if $P(S)=S(1)$, that is, if $S$ is a
Nichols algebra. Now, $\B(V)=\B(X^{-1},\tilde{q})$. Therefore, we
want to show that the relations in $\B(V)=T(V)/\J_V$ hold in $S$.
By hypothesis, we deal only with relations in degree $2$. It is
easily checked that, if $r=b_C\in\J_V^2$ as in
\eqref{eqn:quadraticrelations}, then $r$ is primitive and
satisfies $c(r\ot r)=r\ot r$, where $c$ is the braiding in
$_{G}^{G}\mathcal{YD}$. As $\dim S<\infty$, $r=0$ for any such
$r$. Therefore, there is a surjection $\B(V)\twoheadrightarrow S$,
$P(S)=V=S(1)$ and $S$ is a Nichols algebra.
\end{proof}

\subsection{Lifting of quadratic relations}

Let $X$ be a finite rack, $q$ a 2-cocycle on $X$. Let $H$ be a
pointed Hopf algebra such that its infinitesimal braiding arises
from a principal YD-realization $(\cdot, g, (\chi_i)_{i\in X})$ of
$(X,q)$ over a group $G$. Let $\gr H=\bigoplus_{i>0}
H_i/H_{i-1}=R\sharp\k G$. The projection of Hopf bimodules
$\pi:H_1\to H_1/H_0$ is a morphism of Hopf bimodules over $\k G$.
Let $\sigma$ be a section of Hopf bimodules over $\k G$ and set
$a_i=\sigma(x_i\sharp 1)\in H_1$, then
\begin{equation}\label{nueve}
a_i \ \mbox{is} \ (g_i,1)-\mbox{primitive and} \
g_ia_jg_i^{-1}=q_{ij}a_{i\rhd j}, \ \forall\,i,j\in X.
\end{equation}

Let $\mR'$ be as in Subsect. \ref{subsect:R}; recall that the
space of quadratic relations $\J^2$ for $\B(X,q)$ is generated by
the relations $b_C$ for $C\in\mR'$.  For such $C$, let
\begin{equation}\label{eqn:phiC}
\phi_C = \sum_{h=1}^{n(C)}\eta_h(C) \, X_{i_{h+1}}X_{i_h}
\end{equation}
be a quadratic polynomial in the non-commuting variables $\{X_i :
i\in X\}$, with $\eta_h(C)$ as in \eqref{eqn:quadraticrelations}.
Thus $b_C=\phi_C(\{x_i\}_{i\in X})$. If $(i,j)\in C$, set
$$a_C=\phi_C(\{a_i\}_{i\in X})\in H, \quad  g_C= g_j g_i \in G.$$
Thus the elements $a_C$ are $(g_C,1)$-primitives.

\begin{lem}\label{lem:lambda}
Let $a_i$ be as in \eqref{nueve} and let us suppose that $g$
satisfies:
\begin{equation}\label{eqn:gigj}  g_i\neq
g_jg_k, \quad \forall\,i,j,k\in X.\end{equation} Then there exist
$\lambda_C\in\k$, for all $C\in\mR'$, normalized by
\begin{equation}\label{lambda1}
\lambda_C=0, \ \mbox{if } \ g_C=1, \end{equation} such that:
\begin{align}
\label{lambda2} a_C&=\lambda_C(1-g_C) \text{ in $H$,} \
C\in\mR',\\
\label{lambda1bis} \lambda_C&=q_{k i_{2}}q_{k i_{1}}\lambda_{k\rhd
C},\quad \forall\, k\in X,
\end{align}
if $C=\{i_1,\ldots,i_n\}\in\mR'$ and $k\rhd C=\{k\rhd
i_1,\ldots,k\rhd i_n\}$.

If $X$ is faithful and indecomposable, and $q$ is constant, then the
1-cocycle $(\chi_i)_{i\in X}$ is constant, $\chi_i=\chi$, for all
$i\in X$ and a multiplicative character $\chi$ of $G$, and we have
$\lambda_C=\chi^2(t)\lambda_{t\cdot C}, \ \forall\, t\in G.$
\end{lem}
\begin{proof} Let $C\in \mR'$. As $a_C$ is $(g_C,1)$-primitive,
there are $\lambda_C$ and $\lambda_i\in \k$, for every $i\in X$ such
that $a_C=\lambda_C(1-g_C)+\sum_{i\in X : \, g_i=g_C}\lambda_ia_i$.
Condition \eqref{eqn:gigj} forces $\lambda_i=0$ for every $i\in X$
and thus \eqref{lambda2} follows. Relation \eqref{lambda1bis}
follows by applying $\ad(g_k)$ to both sides of \eqref{lambda2}.
Indeed, let $C'=k\rhd C=\{i_1', \ldots, i_n'\}$ and note first that,
as $i'_1=k\rhd i_1$, $i'_2=k\rhd i_2$, then $i'_l=k\rhd i_l$. Now,
\begin{align*}
  \ad(g_k)(a_C)&=\sum_{h=1}^n\eta_h(C)q_{ki_{h+1}}q_{ki_h}a_{k\rhd
i_{h+1}}a_{k\rhd i_h}\\
&=\sum_{h=1}^n\eta_h(C)q_{ki_{h+1}}q_{ki_h}a_{i'_{h+1}}a_{i'_h}.
\end{align*}
On the other hand,
$$\ad(g_k)(\lambda_C(1-g_C))=\lambda_C(1-g_{C'})=\frac{\lambda_C}{\lambda_{C'}}a_{C'}=\frac{\lambda_C}{\lambda_{C'}}\sum_{h=1}^n\eta_h(C')a_{i'_{h+1}}a_{i'_{h}}.$$
but
\begin{align*}
  \eta_h(C')&=(-1)^{h+1}q_{i'_2i'_1}\ldots q_{i'_hi'_{h-1}}\\
&=(-1)^{h+1}q_{k\rhd i_{2}k\rhd i_{1}}\ldots q_{k\rhd i_{h}k\rhd i_{h-1}}\\
&\overset{\eqref{eqn:condcocycle}}{=}(-1)^{h+1}\prod_{k=1}^{h-1}q_{k i_{k+1}\rhd i_{k}}\prod_{k=1}^{h-1} q_{i_{k+1}i_{k}}\left(\prod_{k=1}^{h-1}q_{ki_{k}}\right)^{-1}\\
&=\eta_{h}(C)\prod_{k=3}^{h+1}q_{k i_k}\left(\prod_{k=1}^{h-1}q_{ki_{k}}\right)^{-1}\\
&=\eta_{h}(C)q_{ki_{h}}q_{ki_{h+1}}q_{ki_{1}}^{-1}q_{ki_{2}}^{-1}.
\end{align*}
Therefore,
$\ad(g_k)(\lambda_C(1-g_C))=\frac{\lambda_C}{\lambda_{C'}}q_{ki_{1}}^{-1}q_{ki_{2}}^{-1}\ad(g_k)(a_C)$.
Now, we have that $\ad(g_k)(a_C)=0\Leftrightarrow
a_C=0\Leftrightarrow g_C=1\Leftrightarrow g_{C'}=1$ and in this case
$\lambda_C=\lambda_{C'}$ by \eqref{lambda1}. If $\ad(g_k)(a_C)\neq
0$, we have \eqref{lambda1bis}. Analogously, the last relation
follows by applying $\ad (H_t)$ to both sides of \eqref{lambda2}:
\begin{align*}
  &\ad (H_t)(\sum_{h=1}^n a_{i_{h+1}}a_{i_h})=\sum_{h=1}^n\chi^2(t)a_{t\cdot i_{h+1}}a_{t\cdot
i_h}=\chi^2(t)a_{t\cdot C},\quad\text{and}\\
&\ad (H_t)(\lambda_C(1-g_C))=\frac{\lambda_C}{\lambda_{t\cdot
C}}\lambda_{t\cdot C}(1-g_{t\cdot
C})=\frac{\lambda_C}{\lambda_{t\cdot C}}a_{t\cdot C},
\end{align*}
thus $\lambda_C=\chi^2(t)\lambda_{t\cdot C}$.
\end{proof}

\begin{cor}\label{cor:lambda}
In the situation above, suppose that $(X,q)$ is $(\mO_2^n,-1)$,
$(\mO_4^4,-1)$ or $(\mO_2^n,\chi)$. Let $a_i$ be as in
\eqref{nueve}. Then $a_C=\lambda_C(1-g_C)$ in $H$, for $C\in\mR'$
and $\lambda_{C}\in\k$ as in Lemma \ref{lem:lambda}.
\end{cor}
\begin{proof}
For $(X,q)=(\mO_2^n,-1)$, see \cite[Lem. 3.4]{AG1}. We follow the
proof there to check condition \eqref{eqn:gigj} in Lemma
\ref{lem:lambda} for $(\mO_4^4,-1)$. In this case, for every
$i,j,k\in X$, $g_i$ acts in the basis $X$ by -1 times a permutation
matrix (since the cocycle is $q\equiv-1$) while $g_kg_l$ acts as a
permutation matrix in the same basis and thus $g_i\neq g_kg_l$.
Consider now $(\mO_2^n,\chi)$. Assume that $i,k,l\in X$ are such
that $g_i=g_kg_l$. Take $j\in X$, then, $i\rhd j=g_i\cdot
j=g_kg_l\cdot j=k\rhd(l\rhd j)$ and therefore $jikl=iklj$,
$\forall\, j\in X$. This would imply $ikl=\id$ in $\s_n$, which is
not possible. Thus, the corollary follows by the previous lemma.
\end{proof}

\begin{rem}\label{rem:K} Note that when dealing with
the racks of transpositions in $\mO_2^n$ or $\mO_4^4$, condition
\eqref{lambda1bis} determine the existence of at most 3 non-zero
scalars $\lambda_C$, say $\lambda_1,\lambda_2, \lambda_3$, where the
subindex is in correspondence with the number of elements in the
class $C$, since $K$ permutes the classes with the same cardinality.

When $G=\s_4$, the infinitesimal braiding of $H$, $V$, is one of
the modules $M(\mO_2^4,\sg\ot\sg)$, $M(\mO_2^4,\sg\ot\eps)$, or
$M(\mO_4^4,\sg\ot\sg)$, by Theorem \ref{teo:ahs}. Given $V$, if
$(X,q)$ are the associated rack and cocycle, the relations from
Corollary \ref{cor:lambda} hold in $H$ for the realization
$(\cdot,\iota, q)$, where $\iota:X\hookrightarrow\s_4$ is the
inclusion, $\cdot:\s_4\times X\to X$ is the action given by
conjugation, and $q$ is either the cocycle $q\equiv-1$ or the
cocycle $q=\chi$ as in Examples \ref{chi}, as appropriate. In this
case, $g_\sigma=H_\sigma$, $\forall\,\sigma\in X$. Therefore, by
the normalization condition \eqref{lambda1}, we will have
$\lambda_1=0$ or $\lambda_2=0$ in Corollary \ref{cor:lambda}
according if we are dealing with $(\mO_2^n,q)$ or $(\mO_4^4,-1)$,
respectively. Furthermore, when the cocycle $q=\chi$, equation
\eqref{lambda1bis} forces $\lambda_2=0$.
\end{rem}

\begin{fed}\label{fed:ql-datum}
A \emph{quadratic lifting datum}, or ql-datum, $\mQ$ consists of
\begin{itemize}
    \item a rack $X$,
    \item a 2-cocycle $q$,
    \item a finite group $G$,
    \item a principal YD-realization $(\cdot, g, (\chi_i)_{i\in X})$ of $(X,q)$ over
    $G$ such that $g$ satisfies \eqref{eqn:gigj},
    \item a collection $(\lambda_C)_{C\in \mR'}$ satisfying
    \eqref{lambda1} and \eqref{lambda1bis}.
    \end{itemize}
\end{fed}

\begin{fed}\label{fed:lifting-of-quadratic}
Given a ql-datum $\mQ$, we define the algebra $\mH(\mQ)$ by
generators $\{a_i, H_t : i\in X, \, t\in G\}$ and relations:
\begin{align}
\label{eqn:hq1} H_e &=1, \quad H_tH_s=H_{ts}, \ t,s\in G;\\
\label{eqn:hq2} H_ta_i &= \chi_i(t)a_{t\cdot i}H_t, \ t\in G, \, i \in X; \\
\label{eqn:hq3} \phi_C(\{a_i\}_{i\in X}) &=
\lambda_C(1-H_{g_jg_i}), \ C\in \mR', \, (i,j)\in C.
\end{align}
Here $i_1=j$, $i_2=i$, $i_{h+2}=i_{h+1}\rhd i_h$ and $\phi_C$ is
as in  \eqref{eqn:phiC}. We shall denote by $a_C$ the left-hand
side of \eqref{eqn:hq3}.
\end{fed}

\begin{rem}\label{rem:pointed}
$\mH(\mQ)$ is a pointed Hopf algebra, setting $\Delta(H_t)=H_t\ot
H_t$,  $\Delta(a_\sigma)=g_\sigma\ot a_\sigma+a_\sigma\ot 1$,
$t\in G$, $\sigma\in X$.  Indeed, it is readily checked that the
comultiplication is well-defined. In this way, $\mH(\mQ)$ is
generated by group-likes and skew-primitives; therefore it is
pointed by \cite[Lem. 5.5.1]{M}.
\end{rem}

\subsection{Pointed Hopf algebras over $\s_n$}

Take $A$ a finite-dimensional pointed Hopf algebra over $\s_4$. As
said in Remark \ref{rem:K}, in this case the ql-datums
\begin{itemize}
    \item $\mQ_n^{-1}[t]=(\s_n,\mO_2^n,-1,\cdot,\iota,\{0,\Lambda,\Gamma\})$,
    \item
    $\mQ_n^{\chi}[\lambda]=(\s_n,\mO_2^n,\chi,\cdot,\iota,\{0,0,\lambda\})$ and
    \item $\mD[t]=(\s_4,\mO_4^4,-1,\cdot,\iota,\{\Lambda,0,\Gamma\})$;
\end{itemize}
for $\Lambda,\Gamma,\lambda\in\k$, $t=(\Lambda,\Gamma)$, are of
particular interest. We will write down the algebras $\hq$ for
these datums in detail. In this case, relations \eqref{lambda2}
for each $C\in\mR'$ with the same cardinality are
$\s_4$-conjugated. Thus it is enough to consider a single relation
for each $C$ with a given number of elements. See Lemma \ref{proj}
for isomorphisms between algebras in the same family.
\begin{fed}\cite[Def. 3.7]{AG1}\label{aene} $\mH(\mQ_n^{-1}[t])$ is the
algebra presented by generators $\{a_i, H_r :
i\in\mO_2^n,r\in\s_n\}$ and relations:
\begin{align*}
&H_e=1, \quad H_rH_s=H_{rs}, \ r,s\in\s_n;\\
&H_j a_i=-a_{jij}H_j, \ i,j\in \mO_2^n; \\
& a_{(12)}^2=0;\\
& a_{(12)} a_{(34)}+a_{(34)} a_{(12)}=\Gamma(1-H_{(12)}H_{(34)});\\
& a_{(12)} a_{(23)}+a_{(23)} a_{(13)}+a_{(13)}
a_{(12)}=\Lambda(1-H_{(12)}H_{(23)}).
\end{align*}
\end{fed}

\begin{fed}\label{achi}
$\mH(\mQ_n^{\chi}[\lambda])$ is the algebra presented by generators
$\{a_i, H_r : i\in\mO_2^n,r\in\s_n\}$ and relations:
\begin{align*}
&H_e=1, \quad H_rH_s=H_{rs}, \ r,s\in\s_n;\\
&H_j a_i=\chi_i(j)a_{jij}H_j, \ i,j\in \mO_2^n; \\
& a_{(12)}^2=0;\\
& a_{(12)} a_{(34)}-a_{(34)} a_{(12)}=0;\\
&a_{(12)}a_{(23)}-a_{(23)}a_{(13)}-a_{(13)}a_{(12)}=\lambda(1-H_{(12)}H_{(23)}).
\end{align*}
\end{fed}

\begin{fed}\label{ace}
$\mH(\mD[t])$ is the algebra generated by elements $\{a_i, H_r :
i\in\mO_4^4,r\in\s_n\}$ and relations:
\begin{align*}
&H_e=1, \quad H_rH_s=H_{rs}, \ r,s\in\s_n;\\
&H_j a_i=-a_{jij}H_j, \ i\in \mO_4^4,j\in\mO_2^4; \\
& a_{(1234)}^2=\Gamma(1-H_{(13)}H_{(24)});\\
& a_{(1234)}a_{(1432)}+a_{(1432)}a_{(1234)}=0;\\
&a_{(1234)}a_{(1243)}+a_{(1243)}a_{(1423)}+a_{(1423)}a_{(1234)}=\Lambda(1-H_{(12)}H_{(13)}).
\end{align*}
\end{fed}

\section{Quadratic algebras}\label{s1}
Finite dimensional Nichols algebras in
$_{\s_4}^{\s_4}\mathcal{YD}$ are all defined by quadratic
relations, see Theorem \ref{teo:ahs}. We need to know the
dimension of the algebras in Definitions \ref{aene}, \ref{achi}
and \ref{ace}, in order to show that they are liftings of these
Nichols algebras. To do this, we first develop a technique on
quadratic algebras to obtain a bound on the dimensions. We follow
\cite{BG} for our exposition. We fix a vector space $W$, and we
let $T=T(W)$ be the tensor algebra. This algebra presents a
natural grading $T=\oplus_{n\geq 0}T^n$, with $T^0=\k$ and
$T^n=W^{\ot n}$ and an induced increasing filtration $F^i$, with
$F^i=\oplus_{j\leq i}T^j$. Let $R\subset W\ot W$ be a subspace and
denote by $J(R)$ the two-sided ideal in $T$ generated by $R$. A
(homogeneous) quadratic algebra $Q(W,R)$ is the quotient
$T(W)/J(R)$. Analogously, for a subspace $P\subset F^2=\k\oplus
W\oplus W\ot W$, we denote by $J(P)$ the two-sided ideal in $P$
generated by $R$. A (nonhomogeneous) quadratic algebra $Q(W,P)$ is
the quotient $T(W)/J(P)$.

Let $A=Q(W,P)$ be a nonhomogeneous quadratic algebra. It inherits an
increasing filtration $A_n$ from $T(W)$. Explicitly,
$A_n=F^n/J(P)\cap F^n$; let $\Gr A=\oplus_{n\geq 0} A_n/A_{n-1}$ be
the associated graded algebra, where $A_{-1}=0$. Consider the
natural projection $\pi: F^2\to W\ot W$ with kernel $F^1$ and set
$R=\pi(P)\subset W\ot W$. Let $B=Q(W,R)$ be the homogeneous
quadratic algebra defined by $R$. Then we have a natural epimorphism
$\rho: B\to \Gr A$. Explicitly, let $\rho': T(W)\to \Gr A$ be
induced by the map $W\hookrightarrow A_1 \twoheadrightarrow
A_1/A_{0}$. Let $x\in R\subset T^2$. Then, there exist $x_0\in\k$,
$x_1\in W$ such that $x-x_1-x_0\in P$ and therefore $ x=x_1+x_0\in
F^2/F^2\cap J(P)=A_2$, thus $\rho'(x)=0\in A_2/A_1$, since
$x_1+x_0\in A_1$. Hence $\rho'$ induces
$\rho:B=T/J(R)\twoheadrightarrow \Gr A$.

\begin{lem}\label{nucleo}
Let $x=x_0-x_1\in P\cap F^1$ be such that $x_0\in\k$ and $x_1\in W$.
Then $\rho(x_1)=0$, so that $\rho$ factors through a morphism of
graded algebras
$$
\rho_x: \widetilde{B}=B/Bx_1B\to \Gr A,
$$
and therefore, $\dim \Gr^nA\leq \dim \widetilde{B}^n$,
$\forall\,n\ge 0$.
\end{lem}
\begin{proof}
Indeed, $x_1+F^1\cap J(P)=x_0\in\k$ and thus $\rho(x_1)=0$. The last
claim follows since from $\rho_x$ is an epimorphism of graded
algebras.
\end{proof}

\begin{pro}\label{pro:quadratic} Let $\mQ$ be a ql-datum and $\mH(\mQ)$ be as in Def. \ref{fed:lifting-of-quadratic}.
Then
\begin{equation*}
  \dim\Gr^n\mH(\mQ)\leq \dim\widehat{\B}_2^n(X,q)\,|G|.
\end{equation*}
In particular, $\dim \mH(\mQ)\leq \dim\widehat{\B}_2(X,q)|G|.$
\end{pro}

For $H=\mH(\mQ_n^{-1}[t])$ or $ \mH(\mQ_n^{\chi}[\lambda])$ and
$n=4,5$, the above proposition implies:
\begin{itemize}   \item $\dim \mH(\mQ_4^{-1}[t]), \dim \mH(\mQ_4^{\chi}[\lambda]), \dim \mH(\mD[t]) \leq 24^3<\infty$,
\item $\dim \mH(\mQ_5^{-1}[t]) , \dim  \mH(\mQ_5^{\chi}[\lambda]) \leq
2^{15}3^55^3<\infty$.
\end{itemize}

\begin{proof} $\hq$ is the nonhomogeneous quadratic algebra
$Q(W,P)$ defined by $W$ and $P$, for $W=\k\{a_i, H_t : i\in X, t\in
G\}$ and $P\subset \k\oplus W\oplus W\ot W$ the subspace generated
by
\begin{align*}
\{H_e-1&, H_t\ot H_s-H_{ts}, H_t\ot a_i-\chi_i(t)a_{t\cdot i}\ot
H_t,\\ &  a_C - \lambda_C1 + \lambda_C H_{g_jg_i}, \ C\in\mR',\,
t,s\in G,\, i\in X\}.
\end{align*} Let $R=\pi(P)$.
Explicitly, $R\subset W\ot W$ is the subspace generated by
\begin{equation*}
\{H_t\ot H_s, \ H_t\ot a_i - \chi_i(t)a_{t\cdot i}\ot H_t, \ a_C,
\ C\in\mR', t,s\in G, i\in X\}.
\end{equation*}
Let $B=Q(W,R)$ be the homogeneous quadratic algebra defined by $W$
and $R$. Let $Y_G $ be the algebra linearly spanned by the set
$\{1, y_t:t\in G\}$, with unit $1$ and multiplication table:
$$ y_ty_s=0, \quad s,t \in G.$$
If $\widehat{\B}_2=\widehat{\B}_2(X,q)$, then $B\cong
\widehat{\B}_2\sharp Y_G$ where $\sharp$ stands for the commutation
relation $(1\sharp y_t)(a_i\sharp 1)=\chi_i(t)a_{t\cdot i}\sharp
y_t$, $(1\sharp 1)(a_i\sharp 1)= a_{i}\sharp 1$, $t\in G$, $i\in X$.  Thus we have an epimorphism
$\rho:\widehat{\B}_2\sharp Y_G\to \Gr \hq$.

Now, note that $P\cap F^1=\k\{H_e-1\}$ and thus, by Lemma
\ref{nucleo}, we have $\rho(H_e)=0$ and an epimorphism
$\rho_e:\widetilde{B}\to \Gr \hq$, with $\widetilde{B}=B/By_eB$. The
commutation relation and the fact that the elements $\{y_t\}_{t\in
G}$ are pairwise orthogonal, give $By_eB=\widehat{\B}_2y_e\subset
B$. This implies $\dim B^n -\dim (\widehat{\B}^n_2 y_e)\geq
\dim\Gr^n \hq$, and, as $\dim B^n=\dim \widehat{\B}^n_2 (|G|+1)$, we
have $\dim \widehat{\B}^n_2 |G|\geq \dim \Gr^n\hq$.
\end{proof}

We now apply Proposition \ref{pro:quadratic} to show that all liftings
of some quadratic Nichols algebras are of the form $\hq$.

\begin{theorem}\label{teo:hq} Let $X$ be a rack, $q$ a 2-cocycle.
Let $H$ be a finite-dimensional pointed Hopf algebra, such that its
infinitesimal braiding is a principal YD-realization $(\cdot, g,
(\chi_i)_{i\in X})$ of $(X,q)$ over $G :=G(H)$ that $g$ satisfies
\eqref{eqn:gigj}. Assume that $\B(X,q)$, $\B(X^{-1},\tilde{q})$ are
quadratic and finite-dimensional. Then there exists a collection
$(\lambda_C)_{C\in\mR'}$ satisfying \eqref{lambda1} and
\eqref{lambda1bis} such that $H\cong \hq$, for the ql-datum
$\mQ=(X,q,G,(\cdot, g, (\chi_i)_{i\in X}),(\lambda_C)_{C\in\mR'})$.
\end{theorem}
\begin{proof}
By Th. \ref{conj}, $\gr (H)\cong \B(X,q)\sharp\k G$. Therefore, we
have $\dim H=\dim\B(X,q)|G|$. On the other hand, by Lemma
\ref{lem:lambda}, there exists a collection
$\{\lambda_C\}_{C\in\mR'}$ and an epimorphism $\hq\twoheadrightarrow
H$ for the ql-datum $\mQ=(X,q,G,(\cdot,g,(\chi_i)_{i\in
X}),\{\lambda_C\}_{C\in\mR'})$. Thus $\dim H\leq \dim\hq$. Now,  by
Prop. \ref{pro:quadratic}, $\dim\hq\leq \dim\B(X,q)|G|$. Therefore,
$H\cong \hq$.
\end{proof}
\section{Representation theory}

We have seen that liftings of some quadratic Nichols algebras are of the form $\hq$
and now we investigate the converse, namely when $\hq$ is actually a lifting.

Let $\mQ=(X,q,G,(\cdot,g,(\chi_i)_{i\in X}),\{\lambda_C\}_{C\in\mR'})$
be a ql-datum; assume that $\B(X,q)$,
$\B(X^{-1},\tilde{q})$ are quadratic and finite-dimensional. Let
$V$ be the corresponding YD-module. By definition the group of
group-likes $G(\hq)$ is a quotient of the group $G$ in the datum.
Therefore, if $\pi:G\twoheadrightarrow G(\hq)$ is the induced
epimorphism, any $\hq$-module $W$ becomes a $\k G$-module by letting
$t\cdot w=\pi(t)\cdot w, \ \forall\, w\in W, t\in G$. Thus
any $\hq$-module is a
direct sum of irreducible $\k G$-modules.
For $i\in X$, set $J_i=\{k\in X : g_i=g_k\}$.

\begin{pro}\label{pro:dimension}
If there is a representation  $\rho:\hq\to\End M$ such that
\begin{enumerate}
\renewcommand{\theenumi}{\roman{enumi}}   \renewcommand{\labelenumi}{(\theenumi)}
  \item\label{item:fiel} $\rho_{|G(\hq)}\circ\pi:G\to\End(M)$, is faithful; \smallbreak
  \item\label{item:dim1}  $\rho(a_{i})\notin\k\rho(G(\hq))$, for all $i\in X$;
and
 \smallbreak
\item\label{item:dim2} the sets $\{\rho(a_j)\}_{j\in J_i}$ are linearly
independent for each $i\in X$;
\end{enumerate}
then $\gr\hq=\B(X,q)\sharp\k G$.
\end{pro}
\begin{proof}
By \eqref{item:fiel},
$G(\hq)\cong G$. Now, $\gr\hq\cong\B(W)\sharp\k G$, for
$\B(W)\in\ydg$ a finite-dimensional Nichols algebra, by Th.
\ref{conj}; and $\dim\B(W)\leq\dim\B(X,q)$ by Prop.
\ref{pro:quadratic}. The map $\B(X,q)\ni x_i\mapsto
\bar{a_i}\in \hq_1/\hq_0$ defines a morphism $\phi:V\to W\in\ydg$.
Conditions \eqref{item:dim1} and \eqref{item:dim2} make $\phi$
injective. Thus, we have a monomorphism
$\B(X,q)\hookrightarrow\B(W)$, see \cite[Cor 3.3]{AS1}; so $V\cong W$.
\end{proof}
\begin{rem} If $i, j\in X$, then  $\rho(a_{i})\notin\k\rho(G)$ implies that
 $\rho(a_{j\rhd i})\notin\k\rho(G)$ by \eqref{eqn:hq2}.
Thus, if $X$ is indecomposable, then \eqref{item:dim1} is equivalent to

\begin{enumerate}
  \item[(ii')] $\exists\,i\in X$ such that
$\rho(a_i)\notin\k\rho(G(\hq))$.
\end{enumerate}
On the other hand, if the realization is faithful (for instance, if $X$ is faithful),
then \eqref{item:dim2} is automatic.
\end{rem}
\subsection{The standard representation of $\s_n$}
Recall that $\s_n$ acts on $\k^n$ permuting the vectors from the
standard basis $\{e_1, \ldots, e_n\}$. The vector $e=e_1+\ldots
+e_n$ is fixed by this action. The orthogonal complement to the
space generated by this vector, with the restricted action from
$\s_n$, is the \emph{standard representation} of $\s_n$, that is
the Specht module $S^{(n-1,1)}$ associated to the partition
$n=n-1+1$. The canonical basis for this space is given by $\{v_1,
\ldots, v_{n-1}\}$ for $v_i=e_i-e_n$. The action on the
corresponding basis is described by:
\begin{align*}
  (i,i+1)\cdot v_j &=v_{(i,i+1)(j)}, \ \mbox{if} \ i<n-1,\\
  (n-1,n)\cdot v_j &=\begin{cases}
    v_j-v_n & \mbox{if} \ j<n,\\
    -v_n & \mbox{if} \ j=n,\\
  \end{cases}
\end{align*}
This representation is faithful. If $n\geq 4$, there is another
$(n-1)$-dimensional representation of $\s_n$, namely the Specht
module $S^{(2,1^{n-2})}\cong S^{(n-1,1)}\ot S^{(1^n)}$ associated to
the partition $n=2+1+\cdots +1$, where $S^{(1^n)}$, the module
associated to the partition $n=1+\cdots+1$, is the sign
representation. For this module, we fix the basis $\{w_i\}$, with
$w_i=v_i\ot z$, $i=1,\ldots,n-1$, if $S^{(1^n)}=\k\{z\}$.

Let now $\mQ$ be one of the ql-data $\mQ_n^{-1}[t]$,
$\mQ_n^{\chi}[\lambda]$, $n\geq 4$, or $\mD[t]$. We will construct:
\begin{itemize}
  \item a $\mH(\mQ_n^{-1}[t])$-module $W(n)$ supported on $S^{(n-1,1)}\oplus S^{(2,1^{n-2})}$;
  \item a $\mH(\mQ_n^{\chi}[\lambda])$-module $U(n)$ supported on $S^{(n-1,1)}$;
  \item a $\mH(\mD[t])$-module $V$ supported on $S^{(3,1)}\oplus
S^{(2,1^{2})}$.
\end{itemize}
Note that for $\mQ=\mQ_n^{\chi}[\lambda]$ a
single irreducible $\k\s_n$ module is enough. This is related to the fact that
these algebras only depend on one parameter $\lambda$. We have used
the computer program \verb"Mathematica"$^\copyright$ to find the
representations. Here we limit ourselves to give a basis $\mB$ of
these modules and write down the base-exchange matrix between $\mB$
and the canonical basis of the $\k\s_4$-module. We also write down
the matrix defining the action of $a_{(12)}$ or $a_{(1234)}$ with
respect to the basis $\mB$. The action of the other elements $a_i$,
$i\in X$ may be deduced from this one through commutation relations
\eqref{eqn:hq2}. Thus, we will have the following proposition, which
provides a Gr\"obner basis-free proof of the analogous result in
\cite[Th. 3.8]{AG1}, for $\mH(\mQ_n^{-1}[t])$.
\begin{pro}\label{cor:dimension}
Let $n\geq 4$, $t\in\k^2, \lambda\in\k$. Let $A_n=\mH(\mQ_n^{-1}[t])$ or
$\mH(\mQ_n^{\chi}[\lambda])$.
\begin{itemize}
\item $G(A_n)\cong\s_n$, $A_n\ncong\k\s_n$.
\item $\dim A_4=24^3$, $\dim A_5 =2^{15}3^55^3$.
\item $G(\mH(\mD[t]))\cong\s_4$, $\mH(\mD[t])\ncong\k\s_4$ and $\dim \mH(\mD[t]) =24^3$.
\end{itemize}
Therefore, for $n=4,5$, the graded Hopf algebra associated to the coradical
filtration of $\mH(\mQ_n^{-1}[t])$ (respectively
$\mH(\mQ_n^{\chi}[\lambda]), \mH(\mD[t])$) is isomorphic to
$\B(\mO_2^n,-1)\sharp \k\s_n$ (respectively $\B(\mO_2^n,\chi)\sharp
\k\s_n$, $\B(\mO_4^4,-1)\sharp \k\s_4)$.
\end{pro}
\begin{proof}
In view of Proposition \ref{pro:dimension}, it follows from Th.
\ref{teo:ahs} and Rem. \ref{rem:O23}, together with Th.
\ref{nichols:cte} and Prop. \ref{nichols:milsch}, using the
representations defined on Propositions \ref{pro:repr An},
\ref{pro:repr Anchi} and \ref{pro:repr D}, as appropriate.
\end{proof}

\subsection{Construction of $\mH(\mQ)$-modules}
Let $t=(\Gamma:\Lambda)\neq(0,0)$.
\begin{pro}\label{pro:repr An}
There is an irreducible $\mH(\mQ_n^{-1}[t])$-module $W(n)$ such that
\begin{itemize}
  \item $W(n)$ has a basis $\{\xi_i,\zeta_i\}_{i=1}^{n-1}$ such that
\begin{equation*}
  \begin{array}{cc}
    a_{12}\xi_1=2\zeta_1, & a_{12}\zeta_1=0,\\
    a_{12}\xi_2=0, & a_{12}\zeta_2=\alpha_n(t) \xi_2,\\
    a_{12}\xi_j=0, & a_{12}\zeta_j=\Gamma \xi_j, j\geq 3.
  \end{array}
\end{equation*}
\item $W(n)\cong_{\s_n}S^{(n-1,1)}\oplus S^{(2,1^{n-2})}$.
%\item $W(n)$ is irreducible.
\end{itemize}
Here, $\alpha_n(t)=2\dfrac{(n-2)\Lambda-(n-3)\Gamma}n$.
\end{pro}
\begin{proof}
Let $n\geq 4$, fix $b=b_n=\frac2{2-n}$, and let $\phi_n:W(n)\to
S^{(n-1,1)}\oplus S^{(2,1^{n-2})}$ be the linear isomorphism
defined, on the basis $\{\xi_i,\zeta_i\}_{i=1}^{n-1}$ and
$\{v_i,w_i\}_{i=1}^{n-1}$ by
\begin{equation*}
\begin{array}{cc}
[\phi_n]=\begin{pmatrix}
\Phi_n & 0  \\
0 & \Phi_n
\end{pmatrix}, & \text{for } [\Phi_n]=\begin{pmatrix}
1 & 1 & 0 &\ldots & 0  \\
-1 & 1 & 0 &\ldots & 0 \\
0 & b & \quad &\quad &\quad \\
\vdots & \vdots & \quad &\id_{n-3} & \quad \\
0 & b & \quad &\quad &\quad \\
\end{pmatrix}.
\end{array}
\end{equation*}
We define a structure of $\k\s_n$-module on $W(n)$ to make $\phi_n$
into an isomorphism of $\k\s_n$-modules. We have thus defined the
action of the elements $H_t$ on this module. For instance, for
$\textbf{0}_{n-1}\in\k^{n-1\times n-1}$ the null matrix; $\rho_{ij}$, for
$2<i<j<n$, the matrix that interchanges the rows $i$ and $j$, we
have
\begin{equation*}
  \begin{array}{cc}
    [H_{(12)}]=\begin{pmatrix}
      \alpha&\textbf{0}_{n-1}\\
\textbf{0}_{n-1}&-\alpha
    \end{pmatrix},&  [H_{(23)}]=\begin{pmatrix}
      \beta&\textbf{0}_{n-1}\\
\textbf{0}_{n-1}&-\beta
    \end{pmatrix},\\
     {[H_{(ij)}]}=\begin{pmatrix}
      \rho_{ij}&\textbf{0}_{n-1}\\
\textbf{0}_{n-1}&-\rho_{ij}
    \end{pmatrix}, & [H_{(n-1 n)}]=\begin{pmatrix}
      \omega&\textbf{0}_{n-1}\\
\textbf{0}_{n-1}&-\omega\end{pmatrix}
  \end{array},
\end{equation*}
where $\alpha_{kl}=\delta_{k,l}\eta_k$, with $\eta_2=-1$,
$\eta_k=1$, $k\neq 2$ and $[\omega]$, $[\beta]$ are, respectively:
\begin{equation*}
\begin{array}{cc}
\begin{pmatrix}
    &&&&0\\
    &\id_{n-1}&&&\vdots\\
        &&&&0\\
    0&0&-1&\ldots&-1
  \end{pmatrix}
,& \begin{pmatrix}\frac12&\frac12-b&-\frac12&0&\ldots&0\\
\frac12&\frac12+b&\frac12&0&\ldots&0\\
-1-b&1+b-2b^2&-b&0&\ldots&0\\
-b&1+b-2b^2&-b& & & \\
\vdots&\vdots&\vdots& &\id_{n-4} &\\
-b&1+b-2b^2&-b& & &
   \end{pmatrix}
   \end{array}
\end{equation*}Thus
$H_{(12)}a_{12}H_{(12)}=-a_{12}=H_{(ij)}a_{12}H_{(ij)}=H_{(n-1n)}a_{12}H_{(n-1n)}$
and the commutativity relations hold. From this matrices, we may compute
\begin{equation*}
  \begin{array}{cc}
    [a_{13}]=\begin{pmatrix}
      \textbf{0}_{n-1} & a_{13}[1]\\
      a_{13}[2]& \textbf{0}_{n-1}
    \end{pmatrix}, &     [a_{23}]=\begin{pmatrix}
      \textbf{0}_{n-1} & a_{23}[1]\\
      a_{23}[2]& \textbf{0}_{n-1}
    \end{pmatrix},
  \end{array}
\end{equation*}
where $a_{13}[1]$ and $a_{23}[1]$ are, respectively, the matrices
\begin{equation*}
\begin{pmatrix}
    \frac\Lambda2&\frac{\Lambda(n-4)+\Gamma(6-2n)}{2(n-2)}&\frac{\Lambda-\Gamma}2&0&\cdots&0\\
        \frac{\Lambda(n-4)+\Gamma(6-2n)}{2n}&\frac{\Lambda(n^2-8n+16)+\Gamma(8n-24)}{2n(n-2)}&\frac{\Lambda(n-4)+\Gamma(6-n)}{2n}&0&\cdots&0\\
            \frac{(n-3)(\Lambda-\Gamma)}{n-2}&\frac{(n-3)(\Lambda(n-4)+\Gamma(6-n))}{(n-2)^2}&\frac{\Lambda(n-3)+\Gamma(4-n)}{n-2}&0&\cdots&0\\
            \frac{\Gamma-\Lambda}{n-2}&\frac{\Lambda(4-n)+\Gamma(n-6)}{(n-2)^2}&\frac{2\Gamma-\Lambda}{n-2}&&&\\
            \vdots&\vdots&\vdots&&\Gamma\id_{n-4}&\\
            \frac{\Gamma-\Lambda}{n-2}&\frac{\Lambda(4-n)+\Gamma(n-6)}{(n-2)^2}&\frac{2\Gamma-\Lambda}{n-2}&&&\\
    \end{pmatrix},
    \end{equation*}
\begin{equation*}
\begin{pmatrix}
    \frac\Lambda2&\frac{\Lambda(4-n)+\Gamma(2n-6)}{2(n-2)}&\frac{\Gamma-\Lambda}{2}&0&\cdots&0\\\\
\frac{\Lambda(4-n)+\Gamma(2n-6)}{2n}&\frac{\Lambda(n^2-8n+16)+\Gamma(8n-24)}{2n(n-2)}&\frac{\Lambda(n-4)+\Gamma(6-n)}{2n}&0&\cdots&0\\\\
\frac{(\Gamma-\Lambda)(n-3)}{n-2}&\frac{(n-3)(\Lambda(n-4)+\Gamma(6-n))}{(n-2)^2}&\frac{\Lambda(n-3)+\Gamma(4-n)}{n-2}&0&\cdots&0\\\\
\frac{\Lambda-\Gamma}{n-2}&\frac{\Lambda(4-n)+\Gamma(n-6)}{(n-2)^2}&\frac{2\Gamma-\Lambda}{n-2}&&&\\\\
            \vdots&\vdots&\vdots&&\Gamma\id_{n-4}&\\\\
\frac{\Lambda-\Gamma}{n-2}&\frac{\Lambda(4-n)+\Gamma(n-6)}{(n-2)^2}&\frac{2\Gamma-\Lambda}{n-2}&&&\\\\
    \end{pmatrix},     \end{equation*}
    and $a_{13}[2], a_{23}[2]$ are such that they both have null $j$-th column for
$j>3$ and have the first three columns as follows:
\begin{equation*}
    \begin{array}{cc}
  a_{13}[2]=\begin{bmatrix}
\frac12&\frac{n}{2(n-2)}&-\frac12\\
        \frac12&\frac{n}{2(n-2)}&-\frac12\\
            \frac{n-3}{2-n}&\frac{n(3-n)}{(n-2)^2}&\frac{3-n}{2-n}\\
            \frac{1}{n-2}&\frac{n}{(n-2)^2}&\frac{1}{2-n}\\
            \vdots&\vdots&\vdots\\
            \frac{1}{n-2}&\frac{n}{(n-2)^2}&\frac{1}{2-n}\\   \end{bmatrix},
      &       a_{23}[2]=\begin{bmatrix}
\frac12&\frac{n}{2(2-n)}&\frac12\\
        -\frac12&\frac{n}{2(n-2)}&-\frac12\\
            \frac{n-3}{n-2}&\frac{n(3-n)}{(n-2)^2}&\frac{n-3}{n-2}\\
\frac{1}{2-n}&\frac{n}{(n-2)^2}&\frac{1}{2-n}\\
            \vdots&\vdots&\vdots\\
\frac{1}{2-n}&\frac{n}{(n-2)^2}&\frac{1}{2-n}\\  \end{bmatrix},
    \end{array}
\end{equation*}
The action of $a_{34}$ differs from the case $n=4$ or $n>4$. If $n=4$,
\begin{equation*}\begin{array}{cccc}
a_{34}\xi_i=2\delta_{i,3}\zeta_i,&a_{34}\zeta_1=\Gamma\xi_1,&a_{34}\zeta_2=(\Lambda-\frac\Gamma2)\xi_2,&a_{34}\zeta_3=0,
\end{array}
\end{equation*}
while if $n>4$ the action is given by
\begin{equation*}
[a_{34}]=\begin{pmatrix}
0&a_{34}[1]\\
a_{34}[2]&0
           \end{pmatrix}
\end{equation*}
where $a_{34}[2]_{ij}=(\delta_{i,3}-\delta_{i,4})(\delta_{j,3}-\delta_{j,4})$ and
\begin{equation*}
 a_{34}[1]=\begin{pmatrix}
0&0&0&0&\ldots&0\\
\frac{8\Lambda+\Gamma(n^2-2n-12)}{n(n-2)}&\frac{3\Gamma-2\Lambda}{n}&\frac{3\Gamma-2\Lambda}{n}&0&\cdots&0\\
\frac{2(n-4)(3\Gamma-2\Lambda)}{(n-2)^2}&\frac{\Lambda(4-n)+\Gamma(n-5)}{2-n}&\frac{\Lambda(4-n)+\Gamma(n-5)}{2-n}&0&\cdots&0\\
\frac{2(n-4)(3\Gamma-2\Lambda)}{(n-2)^2}&\frac{\Lambda(4-n)+\Gamma(n-5)}{2-n}&\frac{\Lambda(4-n)+\Gamma(n-5)}{2-n}&0&\cdots&0\\
\frac{4(2\Lambda-3\Gamma)}{(n-2)^2}&\frac{3\Gamma-2\Lambda}{n-2}&\frac{3\Gamma-2\Lambda}{n-2}&&&\\
\vdots&\vdots&\vdots&&\id_{n-4}&\\
\frac{4(2\Lambda-3\Gamma)}{(n-2)^2}&\frac{3\Gamma-2\Lambda}{n-2}&\frac{3\Gamma-2\Lambda}{n-2}&&&\\
          \end{pmatrix}
\end{equation*}
We can now check the quadratic relations. First, it is clear that
$a_{12}^2=0$. We are left with equations
\begin{equation}\label{eqn:dos}
  a_{12} a_{34}+a_{34}
a_{12}=\Gamma(1-H_{(12)}H_{(34)})
\end{equation}
\begin{equation}\label{eqn:tres}
  a_{12} a_{23}+a_{23}
a_{13}+a_{13} a_{12}=\Lambda(1-H_{(12)}H_{(23)})
\end{equation}
Notice that in both equations, the left and right hand side is zero when computed on $\xi_j$, $\zeta_j$ for $j>3$. The relation is checked on the rest of the generators by a rather straightforward, though tedious, computation. For instance, the
left hand side of \eqref{eqn:tres} applied to $\xi_1$ is
\begin{align*}
&a_{12}(\frac12(\zeta_1-\zeta_2)+\frac{n-3}{n-2}\zeta_3+\frac{1}{2-n}\sum_{j>3}\zeta_j)\\
&\quad+a_{23}(\frac12(\zeta_1+\zeta_2)+\frac{n-3}{2-n}\zeta_3+\frac{1}{n-2}\sum_{j>3}\zeta_j)+2a_{13}\zeta_1=\\
\end{align*}
\begin{align*}&=[-\frac12\alpha_n(t)\zeta_2+\frac{(n-3)\Gamma}{n-2}\zeta_3+\frac{\Gamma}{2-n}\sum_{j>3}\zeta_j]+[\frac{\Lambda+\Gamma(n-3)}{2(n-2)}\xi_1\\
&\quad+\frac{\Lambda(4-n)+\Gamma(n^2-n-6)}{2n(n-2)}\xi_2+\frac{(n-3)(2\Gamma-\Lambda)}{(n-2)^2}\xi_3+\frac{(\Lambda-2\Gamma)}{(n-2)^2}\sum_{j>3}\xi_j\\
&\quad+\frac{n-3}{2-n}(\frac{\Gamma-\Lambda}{2}\xi_1+\frac{\Lambda(n-4)+\Gamma(6-n)}{2n}\xi_2+\frac{\Lambda(n-3)+\Gamma(4-n)}{n-2}\\
&\quad+\frac{2\Gamma-\Lambda}{n-2}\sum_{j>3}\xi_j)+\frac{\Gamma}{n-2}\sum_{j>3}\xi_j]+[2(\frac\Lambda2\xi_1+\frac{\Lambda(n-4)+\Gamma(6-2n)}{2n}\xi_2\\
&\quad+\frac{(n-3)(\Lambda-\Gamma)}{n-2}\xi_3+\frac{\Gamma-\Lambda}{n-2}\sum_{j>3}\xi_j)]\\
&=\frac{3\Lambda}2\xi_1-\frac\Lambda2\xi_2+\frac{(n-3)\Lambda}{n-2}\xi_3+\frac{\Lambda}{2-n}\sum_{j>3}\xi_j.
\end{align*}
And this equal the right hand side, since
$H_{(12)}H_{(23)}\xi_1=\frac12(-\xi_1+\xi_2)+\frac{n-4}{2-n}\xi_3+\frac2{n-2}\sum_{j>3}\xi_j.$
Finally, notice that $a_{12}$ permutes the $\s_n$-modules
$S^{(n-1,1)}$ and $S^{(2,1^{n-2})}$, then $W(n)$ is irreducible.
\end{proof}
\begin{pro}\label{pro:repr Anchi}
There is an irreducible $\mH(\mQ_n^{\chi}[\lambda])$-module $U(n)$ such that
\begin{itemize}
  \item $U(n)\cong_{\s_n}S^{(n-1,1)}$.
  \item $U(n)$ has a basis $\{\xi_i\}_{i=1}^{n-1}$ such that
\begin{equation*}
  \begin{array}{cc}
    a_{12}\xi_1=2\sqrt{-\lambda}\xi_2, & a_{12}\xi_j=0, \ j\geq 3.
  \end{array}
\end{equation*}
%\item $U(n)$ is irreducible.
\end{itemize}
\end{pro}
\begin{proof}
Let $\phi_n:U(n)\to S^{(n-1,1)}$ be the linear isomorphism defined,
on the basis $\{\xi_i\}_{i=1}^{n-1}$ and $\{v_i\}_{i=1}^{n-1}$ by
\begin{equation*}
    [\phi_n]=\left(\begin{matrix}
1 & 1 & 0 & \ldots & 0  \\
1 & -1 & 0 & \ldots & 0  \\
0 & 0 & \quad & \quad & \quad  \\
\vdots & \vdots & \quad & \id_{n-3} &\quad  \\
0 & 0 & \quad & \quad & \quad  \\
\end{matrix}\right).
\end{equation*}
We define a structure of $\k\s_n$-module on $U(n)$ to make $\phi_n$
into an isomorphism of $\k\s_n$-modules. In particular, we have
defined the action of the elements $H_t$ on this module. For
instance, $H_{(12)}\xi_2=-\xi_2$, $H_{(12)}\xi_i=\xi_i$, $i\neq 2$.
The matrices $H_{(ij)}$ with $2<i<j<n$ are permutation matrices that
interchange $a_{ij}$ with $a_{ji}$. When $j=n$, we have that the
$k$th row of $H_{in}$ is the $k$th row of the identity $n-1\times
n-1$ while the $i$th row is $(2\,0\,-1\,\ldots\,-1)$.
% \begin{equation*}
%   \begin{pmatrix}
%            2&0&-1& \ldots &-1
%          \end{pmatrix}.
% \end{equation*}
From these matrices it is readily checked that
$H_{(12)}a_{12}H_{(12)}=-a_{12}$ and $H_{ij}a_{12}H_{ij}=a_{12}$ for
every $i,j\notin\{1,2\}$, and it is immediate that $[a_{12}]^2=0$.
Now, for $\textbf{0}=(0,\ldots,0)\in\k^{n-4}$,
\begin{equation*}
\begin{array}{cc}
 [H_{(13)}]=\begin{pmatrix}
    \frac12&-\frac12&\frac12&\textbf{0}\\
-\frac12&\frac12&\frac12&\textbf{0}\\
1&1&0&\textbf{0}\\
\,^t\textbf{0}&\,^t\textbf{0}&\,^t\textbf{0}&\id_{n-4}
  \end{pmatrix},  [H_{(23)}]=\begin{pmatrix}
    \frac12&\frac12&\frac12&\textbf{0}\\
\frac12&\frac12&-\frac12&\textbf{0}\\
1&-1&0&\textbf{0}\\
\,^t\textbf{0}&\,^t\textbf{0}&\,^t\textbf{0}&\id_{n-4}
  \end{pmatrix}.
\end{array}
\end{equation*}Thus, since
$a_{13}=H_{(23)}a_{12}H_{(23)}$, $a_{23}=-H_{(13)}a_{12}H_{(13)}$,
we have that, on the basis $\{\xi_i\}_{i=1}^{n-1}$, the matrices of
these elements are, respectively,
\begin{equation*}
\begin{pmatrix}
          \frac\Lambda4&-\frac\Lambda4&\frac\Lambda4& 0&\ldots &0\\
           -\frac\Lambda4&\frac\Lambda4&-\frac\Lambda4& 0&\ldots &0\\
           -\frac\Lambda2&\frac\Lambda2&-\frac\Lambda2& 0&\ldots &0\\
           \vdots&\vdots&\vdots&\vdots &\ddots &\vdots\\
           0&0&0& 0&\ldots &0\\
         \end{pmatrix},\quad
\begin{pmatrix}
          \frac\Lambda4&\frac\Lambda4&\frac\Lambda4& 0&\ldots &0\\
           \frac\Lambda4&\frac\Lambda4&\frac\Lambda4& 0&\ldots &0\\
           -\frac\Lambda2&-\frac\Lambda2&-\frac\Lambda2& 0&\ldots &0\\
           \vdots&\vdots&\vdots&\vdots &\ddots &\vdots\\
           0&0&0& 0&\ldots &0\\
         \end{pmatrix}.
\end{equation*}
From this description it is straightforward to check that $[a_{12}a_{23}-a_{23}a_{13}-a_{13}a_{12}]=\lambda[1-H_{(12)}H_{(23)}]$.
% \begin{equation*}
%    [a_{12}a_{23}-a_{23}a_{13}-a_{13}a_{12}]=\lambda[1-H_{(12)}H_{(23)}].
%  \end{equation*}
% \begin{equation*}
%   [a_{12}a_{23}-a_{23}a_{13}-a_{13}a_{12}]=\lambda\begin{pmatrix}
%            \frac12&-\frac12&-\frac12& 0&\ldots &0\\
%            \frac12&\frac32&-\frac12& 0&\ldots &0\\
%            -1&1&1& 0&\ldots &0\\
%            \vdots&\vdots&\vdots&\vdots &\ddots &\vdots\\
%            0&0&0& 0&\ldots &0\\
%          \end{pmatrix}=\lambda[1-H_{(12)}H_{(23)}].
% \end{equation*}
Now, the matrix of $a_{34}$ depends on wether $n=4$ or $n>4$. On
each case:
\begin{equation*}
 \begin{array}{cc}
  [a_{34}]=\begin{pmatrix}
           0&0&0\\
           0&0& 0\\
           -\Lambda&0&0
         \end{pmatrix}, n=4; & [a_{34}]=\begin{pmatrix}
           0&0&0& 0&0&\ldots &0\\
           0 &0&0& 0&0&\ldots &0\\
           0&0&\frac\Lambda2& \frac\Lambda2&0&\ldots &0\\
           0&0&\frac\Lambda2& \frac\Lambda2&0&\ldots &0\\
           0 &0&0& 0&0&\ldots &0\\
           \vdots&\vdots&\vdots&\vdots&\vdots &\ddots &\vdots\\
           0&0&0& 0&0&\ldots &0\\
         \end{pmatrix}, n>4.
         \end{array}
\end{equation*}
and we always have $[a_{12}][a_{34}]=0=[a_{34}][a_{12}]$. $U(n)$ is clearly
irreducible.

\end{proof}

\begin{pro}\label{pro:repr D}
There is an irreducible $\mH(\mD[t])$-module $V$ such that
\begin{itemize}
  \item $V\cong_{\s_4}S^{(3,1)}\oplus S^{(2,1^2)}$.
  \item $V$ has a basis $\{\xi_i,\zeta_i\}_{i=1}^3$ such that
\begin{equation*}
  \begin{array}{cc}
     a_{1234}\xi_1=2\zeta_1, &a_{1234}\zeta_1=\Gamma \xi_1,\\
    a_{1234}\xi_2=2\zeta_2, & a_{1234}\zeta_2=\Gamma \xi_2,\\
    a_{1234}\xi_3=0, & a_{1234}\zeta_3=(\Lambda-\Gamma) \xi_3.
  \end{array}
\end{equation*}
%\item $V$ is irreducible.
\end{itemize}
\end{pro}
\begin{proof}
Let $\phi:V\to S^{(3,1)}\oplus S^{(2,1^2)}$ be the linear
isomorphism defined, on the basis $\{\xi_i,\zeta_i\}_{i=1}^{3}$ and
$\{v_i,w_i\}_{i=1}^{3}$ by
\begin{equation*} [\phi]=\left(\begin{matrix}
1 & 1 & 0 & 0 & 0 & 0 \\
-1 & 1 & 0 & 0 & 0 & 0 \\
0 & -1 & 1 & 0 & 0 & 0 \\
0 & 0 & 0 & 1 & 1 & 0 \\
0 & 0 & 0 & -1 & 1 & 0 \\
0 & 0 & 0 & 0 & -1 & 1 \\
\end{matrix}\right).\end{equation*}
We define a structure of $\k\s_4$-module on $V$ to make $\phi$ into
an isomorphism of $\k\s_4$-modules. We have thus defined the action
of the elements $H_t$ on this module. From the definition of the
action for $a_{(1234)}$ it is easy to describe the action of
$a_{(1432)}, a_{(1243)}, a_{(1423)}$ and to check that it defines in
fact an $\mH(\mD[t])$-module. Since $a_{(1234)}$ permutes
$S^{(3,1)}$ and $S^{(2,1^{2})}$, $V$ is irreducible.
\end{proof}
\section{Pointed Hopf algebras over $\s_4$}
In this section we prove the Main Theorem. We first need the following lemma, which
establishes the isomorphism classes between the algebras belonging
to one of the families defined on Definitions
\ref{aene}--\ref{ace}.

\begin{lem}\label{proj}
Let $t, t'\in\k^2$. Then $\mH(\mQ_n^{-1}[t])\cong
\mH(\mQ_n^{-1}[t'])$ if and only if $t\neq 0$ and
$t=t'\in\mathbb{P}_\k^1$ or if $t=t'=(0,0)$, and the same holds for
$\mH(\mD[t])$. Finally, $\mH(\mQ_n^{\chi}[\lambda])\cong
\mH(\mQ_n^{\chi}[1]),\,\forall\,\lambda\in\k^*$ and
$\mH(\mQ_n^{\chi}[1])\ncong \mH(\mQ_n^{\chi}[0])$.
\end{lem}
\begin{proof}
Let $\mQ=\mQ_n^{-1}[t]$, $\mQ_n^{\chi}[\lambda]$ or $\mD[t]$;
$\mQ'=\mQ_n^{-1}[t']$, $\mQ_n^{\chi}[\lambda']$ or $\mD[t']$. Let
$H=\hq$, $H=\mH(\mQ')$. An isomorphism $\phi:H\to H'$ induces
$\phi_{|\s_n}\in\aut(\s_n)$ for $n=4,5$ as appropriate, and thus we
have a permutation $\varphi:X\to X$ such that
$\phi(H_i)=H_{\varphi(i)}$ and $\phi(a_i)$ is
$(H_{\varphi(i)},1)$-primitive, for $i\in X$. Furthermore, $\phi$
restricts to $\phi_1:H_1\to H'_1$ between the second factors from
the coradical filtration, and thus $\phi(a_i)=\eta
a'_{\varphi(i)}+\mu(1-H_{\varphi(i)})$, $\eta,\mu\in\k$. When
$X=\mO_2^n$, relation $a_i^2=0$ makes $\mu=0$, and for $X=\mO_4^4$
relation $a_i^2=\Gamma(1-H_i^2)$ makes $\mu=0$ and
$\Gamma=\eta^2\Gamma'$. Therefore, $\phi(a_i)=\eta a_{\varphi(i)}$.
Now, taking into account the rest of the quadratic relations, we
have that, for $\mQ_n^{-1}[t]$ as well as for $\mD[t]$, $t=\eta^2
t'$; while for $\mQ_n^{-1}[\lambda]$ we get
$\lambda=\eta^2\lambda'$. Thus the result follows.
\end{proof}
\begin{proofmt}
In Remark \ref{rem:pointed} we have already seen in that the
algebras listed are pointed Hopf algebras. We have computed their
dimension in Proposition \ref{cor:dimension}. Reciprocally, let $H$
be a finite-dimensional pointed Hopf algebra with $G(H)\cong \s_4$.
By Th. \ref{teo:ahs}, the infinitesimal braiding of $H$ is
isomorphic either to $M(\mO_2^4,\sg\ot\sg)$, to
$M(\mO_2^4,\sg\ot\eps)$ or to $M(\mO_4^4,\rho_{\_})$. Since all
these modules are self-dual, $H$ is isomorphic to one, and only one,
of the Hopf algebras in the list by Th. \ref{teo:hq}, Rem.
\ref{rem:K} and Lemma \ref{proj}.\qed
\end{proofmt}

\begin{rem} The classification of finite-dimensional pointed Hopf algebras with $G(H)\cong \s_3$,
done in \cite{AHS} using \cite{AG1}, can be alternatively proved in this way.
\end{rem}

We end this paper by describing in the following corollary a
sub-class of finite-dimensional pointed Hopf algebras over $\s_5$.
It follows in the same fashion as the Main Theorem.
\begin{cor}
Let $A$ be a finite-dimensional pointed Hopf algebra with
$G(A)\cong\s_5$ and let $M\in\,_{\s_5}^{\s_5}\mathcal{YD}$ be its
infinitesimal braiding.
\begin{itemize}
\item If $M\cong M(\mO_2^5,\sg\ot\sg)$, then $A\cong\mH(\mQ_5^{-1}[t])$, for exactly one
$t\in\mathbb{P}_\k^1\cup\{(0,0)\}$.
\item If $M\cong M(\mO_2^5,\sg\ot\eps)$, then $A\cong\mH(\mQ_5^{\chi}[\lambda])$, for one $\lambda\in\{0,1\}$.\qed
\end{itemize}
\end{cor}

\subsection*{Acknowledgements} We thank Mat\'ias Gra\~na for his
assistance on the proof of Proposition \ref{nichols:milsch}. A.G.I.
specially thanks Nicol\'as Andruskiewitsch for his patient guidance
and careful reading of this work, which lead to a significant
enhancement on its presentation, accuracy and extent.

\end{document}